\documentclass[12pt]{amsart}

\usepackage{amsmath}
\usepackage{amsfonts}
\usepackage{amssymb}
\usepackage{amscd}
\usepackage{graphicx}
\usepackage[abbrev,alphabetic]{amsrefs}
\RequirePackage[dvipsnames,usenames]{color}
\usepackage{soul,xcolor}
\setstcolor{red}
\usepackage{stmaryrd}
\usepackage{mathtools}
\usepackage{booktabs}
\usepackage{multirow}
\newtagform{tiny}{\tiny(}{)}
\usepackage{bm}

\usepackage{mathtools}
\usepackage{hyperref}
\usepackage[margin=1.25in]{geometry}

\usepackage{cleveref}

\usepackage{amsthm}
\usepackage{comment}
\usepackage[all,cmtip]{xy}
\usepackage{tikz-cd}
\usetikzlibrary{cd}

\usepackage[all]{xy}

\makeatletter
\def\@tocline#1#2#3#4#5#6#7{\relax
  \ifnum #1>\c@tocdepth % then omit
  \else
    \par \addpenalty\@secpenalty\addvspace{#2}%
    \begingroup \hyphenpenalty\@M
    \@ifempty{#4}{%
      \@tempdima\csname r@tocindent\number#1\endcsname\relax
    }{%
      \@tempdima#4\relax
    }%
    \parindent\z@ \leftskip#3\relax \advance\leftskip\@tempdima\relax
    \rightskip\@pnumwidth plus4em \parfillskip-\@pnumwidth
    #5\leavevmode\hskip-\@tempdima
      \ifcase #1
       \or\or \hskip 1em \or \hskip 2em \else \hskip 3em \fi%
      #6\nobreak\relax
    \hfill\hbox to\@pnumwidth{\@tocpagenum{#7}}\par% <---- \dotfill -> \hfill
    \nobreak
    \endgroup
  \fi}
\makeatother

\newsavebox{\pullback}
\sbox\pullback{%
\begin{tikzpicture}%
\draw (0,0) -- (1ex,0ex);%
\draw (1ex,0ex) -- (1ex,1ex);%
\end{tikzpicture}}

\newsavebox{\pullbackdl}
\sbox\pullbackdl{%
\begin{tikzpicture}%
\draw (-1ex,0ex) -- (0ex,0ex);%
\draw (0ex,-1ex) -- (0ex,0ex);%
\end{tikzpicture}}

\newsavebox{\pushoutdr}
\sbox\pushoutdr{%
\begin{tikzpicture}%
\draw (-1ex,-1ex) -- (-1ex,0ex);%
\draw (-1ex,0ex) -- (0ex,0ex);%
\end{tikzpicture}}

\renewcommand{\mod}{\ \textrm{mod}\ }

\newcommand{\Z}{\mathbb{Z}}

\newcommand{\m}{\mathfrak{m}}

\newcommand{\wt}{\widetilde}

\DeclareMathOperator{\Ker}{Ker}

\theoremstyle{plain}
\newtheorem{theorem}{Theorem}[section]

\newtheorem{proposition}[theorem]{Proposition}
\newtheorem{lemma}[theorem]{Lemma}
\newtheorem{corollary}[theorem]{Corollary}

\newtheorem{claim}[theorem]{Claim}
\newtheorem*{claim*}{Claim}

\newtheorem{theoremA}{Theorem}

\theoremstyle{definition}
\newtheorem{definition}[theorem]{Definition}

\newtheorem{example}[theorem]{Example}
\newtheorem{notation}[theorem]{Notation}

\newtheorem*{setup*}{Setup}
\newtheorem{question}[theorem]{Question}

\theoremstyle{remark}
\newtheorem{remark}[theorem]{Remark}
\newtheorem*{ackn}{Acknowledgements}

\newenvironment{claimproof}[0]
  {%
   \paragraph{\it Proof.}%
  }
  {%
    \hfill$\blacksquare$%
  }

\crefname{theorem}{Theorem}{Theorems}
\crefname{proposition}{Proposition}{Propositions}
\crefname{definition}{Definition}{Definitions}
\crefname{example}{Example}{Examples}
\crefname{lemma}{Lemma}{Lemmas}
\crefname{corollary}{Corollary}{Corollaries}
\crefname{conjecture}{Conjecture}{Conjectures}
\crefname{claim}{Claim}{Claims}
\crefname{notation}{Notation}{Notations}



\numberwithin{equation}{theorem}

\title[Tight closure of ideals on Witt rings]
{Tight closure of ideals on Witt rings}%: local theory}

\author{Shou Yoshikawa}
\address{Tokyo Institute of Technology, Tokyo 152-8551, Japan}
\email{yoshikawa.s.al@m.titech.ac.jp}

\begin{document}

\begin{abstract}
In this paper, we introduce the notions of tight closure of ideals on Witt rings and quasi-tightly closedness of system of parameters.
By using the notions, we obtain a characterization of quasi-$F$-rationality introduced in \cite{KTTWYY3}.
Furthermore, we study the relationship between the closure operator and integrally closure.
\end{abstract}

\subjclass[2010]{}   

\maketitle

\setcounter{tocdepth}{2}

\tableofcontents
\section{Introduction}
One of the fundamental topics of commutative algebra in positive characteristic is the study of Frobenius splitting and Frobenius singularities.
In the theory of Frobenius singularities, the notion of tight closure plays an essential role.
For an ideal $I$ of an $F$-finite Noetherian reduced ring $R$ in characteristic $p>0$, we say that $a \in I^*$, the \emph{tight closure} of $I$, if there exists a non-zero divisor $c \in R$ such that $ca^{p^e} \in I^{[p^e]}$ for every $e \geq 0$ and $I$ is \emph{tightly closed} if we have $I=I^*$.
Furthermore, we say that $R$ is \emph{$F$-regular} if every ideal of $R$ is tightly closed.
\cite{HH90} proved that regularity implies $F$-regularity and if $I$ is generated by $n$ elements, then the integral closure of $I^n$ is contained in $I^*$.
Combining such two results, if $R$ is regular and $I$ is generated by $n$ elements, then the integral closure of $I^n$ is contained in $I$ (the Brian\c{s}on-Skoda theorem).
Moreover, we say that $R$ is \emph{$F$-rational} if every parameter ideal of $R$ is tightly closed (cf.~ \cite{WF87}, \cite{Smith94}).

Recently,  \cite{Yob23} introduce a new notion, called \emph{quasi-$F$-splitting}, motivated by the theory of crystalline cohomology. 
Furthermore, in \cite{KTTWYY1} and \cite{KTTWYY2}, 
we study quasi-$F$-splitting in the context of birational geometry.
We refer to the introductions of \cite{Yob23} and \cite{KTTWYY1} for more information on quasi-$F$-splittings, and to \cite{KTY22} for Fedder type criterion for quasi-$F$-splittings.
Moreover, \cite{TWY}, \cite{KTTWYY3} introduce \emph{quasi-$F$-regularity} and \emph{quasi-$F$-rationality}
which are generalizations of some characterization of $F$-regularity and $F$-rationality, respectively.
We recall the definition of quasi-$F$-rationality.
\begin{definition}
Let $(R,\m)$ be a $d$-dimensional $F$-finite Noetherian local domain of characteristic $p>0$ and $h \geq 1$ an integer.
\begin{itemize}
    \item We define the submodule $\wt{0^*_h}$ of $H^d_\m(W_h(R))$ as follows.
    We say that $\alpha \in \wt{0^*_h}$ if there exists a non-zero divisor $c \in R$ such that $[c]F^e(\alpha)=0$ for every integer $e \geq 0$, where $W_h(R)$ is the Witt vector of $R$ of length $h$.
    \item $R$ is $h$-quasi-$F$-rational if $R$ is Cohen-Macaulay and the image of $\wt{0^*_h}$ by the map $H^d_\m(W_h(R)) \to H^d_\m(R)$ is zero.
\end{itemize}
\end{definition}
\noindent
It is known that $1$-quasi-$F$-rationality is equivalent to $F$-rationality.

In this paper, we first introduce tight closure of ideals on Witt rings as follows.
Let $I_n$ be an ideal of $W_n(R)$,  then we say that $\alpha \in I_n^*$, the \emph{tight closure} of $I_n$, if there exists a non-zero divisor $c \in R$ such that $[c]F^e(\alpha) \in I_n^{[p^e]}$, where $I_n^{[p^e]}$ is an ideal of $W_n(R)$ generated by $F^e(I_n)$.
For $n=1$, it coincides with the original tight closure.
Moreover, for an ideal $I_n \subseteq W_n(R)$, we say that $I_n$ is \emph{tightly closed} if $I_n^*=I_n$.
Furthermore, we study the relationship between plus closure and tight closure for parameter ideals of $W_n(R)$ (\cref{prop:tight-vs-plus}), which is a generalization of \cite{Smith94}*{Proposition~5.1}.
Next, we introduce quasi-tightly closedness of a part of system of parameters as follows:
\begin{definition}[\cref{defn:q-t.c}]\label{defn:Intro-q-t.c}
Let $(R,\m)$ be an $F$-finite Noetherian local domain of characteristic $p>0$ and $h \geq 1$  an integer.
Let $x_1,\ldots,x_r$ be a part of a system of parameters of $R$.
We set $I_h^{\bm{v}}:=([x_1^{v_1}],\ldots,[x_r^{v_r}]) \subseteq W_h(R)$ for every integers $v_1,\ldots,v_r \geq 1$.
We say that $x_1,\ldots,x_r$ is \emph{$h$-quasi-tightly closed} if for every integers $v_1,\ldots,v_r \geq 1$, then we have
\[
(I_h^{\bm{v}})^*R=(x_1^{v_1},\ldots,x_d^{v_d}).
\]
We note that the ideal $(I_n^{\bm v})^*R$ depends on the choice of $x_1,\ldots,x_r$ (cf.~\cref{ex} (2)), thus we do not say "$I$ is quasi-tightly closed''.
\end{definition}
\noindent
Furthermore, we obtain a characterization of quasi-$F$-rationality and the condition (2) in \cref{intro:thm:new-old-h} is an analog of the definition of $F$-rationality.
\begin{theoremA}[\cref{thm:new-old-h}]\label{intro:thm:new-old-h}
Let $(R,\m)$ be a $d$-dimensional $F$-finite Noetherian local domain of characteristic $p>0$ and $h \geq 1$  an integer.
Then the following are equivalent to each other.
\begin{enumerate}
    \item  $R$ is $h$-quasi-$F$-rational, 
    \item for every system of parameters $x_1,\ldots,x_d$  is $h$-quasi-tightly closed.
    \item there exists a system of parameters $x_1,\ldots,x_d$ of $R$ such that $x_1,\ldots,x_d$ is $h$-quasi-tightly closed.
\end{enumerate}
\end{theoremA}
\noindent
The difficulty of the proof of \cref{intro:thm:new-old-h} is to prove Cohen-Macaulayness by assuming the condition (3).
One of the important problem in the theory of quasi-Frobenius singularities is that quasi-$F$-regularity implies Cohen-Macaulayness or not.
We expect that \cref{intro:thm:new-old-h} plays an essential role to study the problem.
By \cref{prop:gene t-c vs i-c}, the Cohen-Macaulayness of a quasi-$F$-regular ring $R$ is reduced to the following question.
\begin{question}
Let $(R,\m)$ be a $d$-dimensional $F$-finite Noetherian local domain of characteristic $p>0$ and $h \geq 1$ an integer.
Let $x_1,\ldots,x_d$ be a system of parameters and $I_h:=([x_1],\ldots,[x_d]) \subseteq W_h(R)$.
Is the map 
\[\omega_R(-K_R) \otimes W_h(R)/I_h \to W_h\omega_R(-K_R) \otimes W_h(R)/I_h
\]
injective, where the map is induced by the trace map $\omega_R \to W_h\omega_R$ (cf.~\cite{KTTWYY3}*{Section~2.4})?
\end{question}

Finally, we obtain the following result, which is a relationship between integral closure and a new closure operator appearing in \cref{defn:Intro-q-t.c}.
Theorem \ref{thm:intro-gene t-c vs i-c} were inspired by a discussion with Jakub Witaszek.
\begin{theoremA}[\cref{prop:int-cl vs tight-cl}]\label{thm:intro-gene t-c vs i-c}
Let $R$ be an $F$-finite Noetherian ring of characteristic $p >0$ and $n \geq 1$ an integer.
Let $f_1,\ldots,f_r \in R$, $I_n:=([f_1],\ldots,[f_r]) \subseteq W_n(R)$ and $I:=(f_1,\ldots,f_r)$.
Then we have $\overline{I^{(\lambda+r-1)}} \subseteq (I_n^{\lambda})^*R$ for every integer $ \lambda \geq 1$.
\end{theoremA}
\noindent
Combining \cref{intro:thm:new-old-h} and \cref{thm:intro-gene t-c vs i-c}, we obtain a Brian\c{s}on-Skoda type theorem for parameter ideals on quasi-$F$-rational rings.
Since quasi-$F$-rationality implies pseudo-rationality (\cite{KTTWYY3}*{Theorem~3.44}), the result is still known by \cite{LT}*{Corollary~2.2}.

\begin{ackn}
The author wishs to express our gratitude to Tatsuro Kawakami, Hiromu Tanaka, Teppei Takamatsu, Fuetaro Yobuko and Jakub Witaszek for valuable discussion.
The author also grateful to Shunsuke Takagi, Kenta Sato for helpful comments.
Finally, the author would like to express their sincere gratitude to the referee for carefully reviewing the manuscript and providing valuable comments.
The author was supported by JSPS KAKENHI Grant number JP20J11886.
\end{ackn}

\section{Big Cohen-Macaulayness of $W_n(R^+)$}

\begin{definition}
Let $R$ be a Noetherian domain with fractional field $K$.
An \emph{absolute integral closure} of $R$ is the integral closure of $R$ in some algebraic closure $\overline{K}$ of $K$.
\end{definition}

\begin{proposition}\label{prop:big-CM}
Let $(R,\m)$ be an $F$-finite Noetherian local domain of positive characteristic.
Let $R^+$ be an absolute integral closure of $R$.
Then $W_n(R^+)$ is a big Cohen-Macaulay $W_n(R)$-algebra for every integer $n \geq 1$, that is, every system of parameters of $W_n(R)$ is a regular sequence on $W_n(R^+)$.
\end{proposition}

\begin{proof}
We note that $W_n(R)$ is a local ring and the unique maximal ideal of $W_n(R)$ is 
\[
\{(a_0,a_1,\ldots,a_{n-1}) \in W_n(R) \mid a_0 \in \m\}. 
\]
Let us prove the assertion by induction on $n$.
For $n=1$, it is due to \cite{HL}*{Corollary 2.3}.
We note that by \cite{Gabber}*{Remark 13.6}, $R$ is a surjective image of some regular ring.
We assume $n \geq 2$.
We take a system of parameters $\alpha_1,\ldots,\alpha_d$ of $W_n(R)$.
Then $F(R(\alpha_1)),\ldots,F(R(\alpha_d))$ and $R^{n-1}(\alpha_1),\ldots,R^{n-1}(\alpha_d)$ are  system of parameters of $W_{n-1}(R)$ and $R$, respectively.
By the induction hypothesis, $F(R(\alpha_1)),\ldots,F(R(\alpha_d))$ and $R^{n-1}(\alpha_1),\ldots,R^{n-1}(\alpha_d)$ are regular sequences of $W_{n-1}(R^+)$ and $R^+$, respectively.
Consider the following commutative diagram in which each horizontal sequence is exact:
\[
\begin{tikzcd}
    0 \arrow[r] & F_*W_{n-1}(R^+)  \arrow[r,"V"] \arrow[d,"\cdot F_*F(\alpha_1) "] & W_n(R^+) \arrow[r,"R^{n-1}"] \arrow[d,"\cdot \alpha_1"] & R^+ \arrow[r] \arrow[d,"\cdot \alpha_1"] & 0 \\
    0 \arrow[r] & F_*W_{n-1}(R^+) \arrow[r,"V"] & W_n(R^+) \arrow[r,"R^{n-1}"] & R^+ \arrow[r] & 0,
\end{tikzcd}
\]
where the left vertical map and the right vertical map are injective.
By the snake lemma, we have $\alpha_1$ is a regular element of $W_n(R^+)$ and the exact sequence
\[
0 \to F_*W_{n-1}(R^+)/F(\alpha_1)W_{n-1}(R^+) \to W_n(R^+)/\alpha_1W_n(R^+) \to R^+/\alpha_1R^+ \to 0.
\]
By inductive argument, we have $\alpha_1,\ldots,\alpha_r$ is a regular element on $W_n(R^+)$ and the exact sequence
\[
0 \to F_*W_{n-1}(R^+)/I_r^{[p]}W_{n-1}(R^+) \to W_n(R^+)/I_rW_n(R^+) \to R^+/I_rR^+ \to 0
\]
for every integers $1 \leq r \leq d$, where $I_r:=(\alpha_1,\ldots,\alpha_r) \subseteq W_n(R)$ and $I_r^{[p]}:=(F(\alpha_1),\ldots,F(\alpha_r)) \subseteq W_n(R)$.
In particular, we have $\alpha_1,\ldots,\alpha_d$ is a regular sequence on $W_n(R^+)$.
\end{proof}

\section{An analog of tight closure of ideals via Witt rings}
In this section, we introduce the notion of tight closure of ideals on Witt rings.
Furthermore, by using it, we introduce a new closure operator, which is called quasi-tight closure and related to quasi-$F$-singularities (cf.~\cref{thm:new-old-h}).

\begin{definition}
Let $R$ be an $F$-finite Noetherian domain of characteristic $p>0$.
Let $I$ be an ideal of $W_n(R)$.
\begin{itemize}
    \item We denote the ideal of $W_n(R)$ generated by $F^e(I)$ by $I^{[p^e]}$.
    \item An ideal $I^*$, the \emph{tight closure} of $I$, is defined as follows.
    Let $\alpha \in W_n(R)$.
    Then $\alpha$ is contained in $I^*$ if and only if there exists $c \in R^\circ$ such that for every integer $e \geq 0$,
    \[
    [c]F^e(\alpha) \in I^{[p^e]}.
    \]
    \item We define the ideal $I^+$ of $W_n(R)$ as $I^+:=IW_n(R^+) \cap W_n(R)$.
\end{itemize}
\end{definition}
 
\begin{notation}\label{n-local}
Let $(R, \m)$ be a Noetherian $F$-finite local domain of characteristic $p>0$. 
Set $d := \dim R$.  
Let $R^\circ$ be the set of non-zero divisors of $R$ and $R^+$ an absolute integral closure of $R$.
\end{notation}

\subsection{Tight closure versus plus closure on Witt rings}

In this subsection, we study a relationship between tight closure and plus closure on Witt rings, which is a generalization of \cite{Smith94}*{Proposition~5.1}.
For a domain $R$ and an ideal $I_n \subseteq W_n(R)$, the \emph{plus-closure} $I_n^+$ of $I_n$ is defined by
\[
I_n^+:=I_n W_n(R^+) \cap W_n(R), 
\]
where $R^+$ denotes an absolute integrally closure of $R$.

\begin{lemma}\label{regular sequence case}
We use \cref{n-local}.
Let $x_1,\ldots, x_r$ be a part of a system of parameters of $R$.
Let $I_n:=([x_1],\ldots, [x_r]) \subseteq W_n(R)$  for integers $n \geq 1$.
\begin{enumerate}
\item If $n \geq 2$, then we have an exact sequence
\begin{equation*}
0 \to F_*W_{n-1}(R^+)/I_{n-1}^{[p^{e+1}]}W_{n-1}(R^+) \to W_n(R^+)/I_n^{[p^e]}W_n(R^+) \to R^+/I_1^{[p^e]}R^+ \to 0
\end{equation*}
for every integer $e \geq 0$.
\item Let $c \in R^\circ$ be a test element.
For every $\alpha \in I_n^+$, $e \geq 0$ and $n \geq 1$, we have
\[
[c^2]F^e(\alpha) \in I_{n}^{[p^e]}.
\]
\end{enumerate}
\end{lemma}

\begin{proof}
Let us show (1) by induction on $r$.
Let $J_n:=([x_1],\ldots,[x_{r-1}])W_n(R^+) \subseteq W_n(R^+)$, where if $r=1$, then $J_n=0$.
By induction hypothesis, we obtain the following commutative diagram in which each horizontal sequence is exact:
\begin{equation} 
    \begin{tikzcd}
    0 \arrow{r} & F_*W_{n-1}(R^+)/J_{n-1}^{[p^{e+1}]} \arrow{r}{V} \arrow{d}{\cdot [x_r^{p^{e+1}}])} & W_n(R^+)/J_n^{[p^e]} \arrow{r} \arrow{d}{\cdot [x_r^{p^e}]} & R^+/J_1^{[p^e]} \arrow{r}{V} \arrow{d}{\cdot x_r^{p^e}} & 0 \\
    0 \arrow{r} & F_*W_{n-1}(R^+)/J_{n-1}^{[p^{e+1}]} \arrow{r} & W_n(R^+)/J_n^{[p^e]} \arrow{r} & R^+/J_1^{[p^e]} \arrow{r} & 0.
    \end{tikzcd}
    \end{equation}
Since $x_1^{p^e},\ldots,x_r^{p^e}$ is a regular sequence on $R^+$ by Proposition \ref{prop:big-CM}, the right vertical map is injective.
Thus, we obtain the following exact sequence by the snake lemma
\[
0 \to F_*W_{n-1}(R^+)/I_{n-1}^{[p^{e+1}]}W_{n-1}(R^+) \to W_n(R^+)/I_n^{[p^e]}W_n(R^+) \to R^+/I_1^{[p^e]}R^+ \to 0,
\]
as desired.

Next, we prove (2) by induction on $n$.
In the case of $n=1$, it is known that $I_1R^+ \cap R \subseteq I_1^*$ by \cite{Smith94}*{Proposition~2.14}.
Thus, we have
\[
c(I_1^{[p^e]})^+ \subseteq c(I_1^{[p^e]})^* \subseteq I_1^{[p^e]},
\]
as desired.
We assume $n \geq 2$ and take $\alpha \in W_n(R) \cap I_nW_n(R^+)$.
We fix an integer $e \geq 0$.
By the result in the case of $n=1$, we have $R^{n-1}([c]F^e(\alpha)) \in I_1^{[p^e]}$.
By the exact sequence
\[
F_*W_{n-1}(R)/I_{n-1}^{[p^{e+1}]} \to W_n(R)/I_n^{[p^e]} \to R/I_1^{[p^e]} \to 0,
\]
there exists $\beta_e \in W_{n-1}(R)$ such that
\[
[c]F^e(\alpha) \equiv V\beta_e \mod I_n^{[p^e]}.
\]
Thus, we have $V\beta_e \in I_n^{[p^e]}W_n(R^+)$, so $\beta_e \in I_{n-1}^{[p^{e+1}]}W_n(R^+)$ by (1).
By the induction hypothesis, we have
\[
[c^2]\beta_e \in I_{n-1}^{[p^{e+1}]}.
\]
Thus, we have
\[
[c]V\beta_e=V([c^p]\beta_e) \in I_n^{[p^e]},
\]
and in particular, we obtain $[c^2]F^e(\alpha) \in I_n^{[p^e]}$.
\end{proof}

\begin{lemma}\label{colon capturing}
We use \cref{n-local}.
Let $x_1,\ldots, x_r$ be a part of a system of parameters of $R$.
Let $I_n:=([x_1],\ldots, [x_r]) \subseteq W_n(R)$  for integers $n \geq 1$.
Then we have the following;
\begin{enumerate}
    \item We have $I_n^+\subseteq I_n^*$ for integers $n \geq 1$.
    \item We have $V^{-1}(I_n^*)=F_*(I_{n-1}^{[p]})^*$ for integers $n \geq 2$.
    \item If $x_1,\ldots,x_r,x_{r+1}$ be a part of a system of parameters, then 
    \[
    (I_n \colon [x_{r+1}]) \subseteq (I_n^* \colon [x_{r+1}])=I_n^*
    \]
    for integers $n \geq 1$.
\end{enumerate}
\end{lemma}

\begin{proof}
By Lemma~\ref{regular sequence case}~(2), we obtain the assertion (1).

Next, let us show the assertion (2).
Take $\beta \in (I_{n-1}^{[p]})^*$, then there exists $c' \in R^\circ$ such that $[c']F^e(\beta) \in I_{n-1}^{[p^{e+1}]}$ for every integer $e \geq 0$.
Then we have $[c'^p]F^eV\beta \in I_n^{[p^e]}$, thus $V\beta \in I_n^*$.
Next, we take $\alpha \in V^{-1}(I_n^*)$.
Then there exists $c' \in R^\circ$ such that 
\[
V([c']F^e(\alpha))=[c'^p]VF^e(\alpha) \in I_n^{[p^e]}
\]
for every $e \geq 0$.
Since we have
\[
V^{-1}(I_n^{[p^e]}) \subseteq F_*W_{n-1}(R) \cap V^{-1}(I_n^{[p^e]}W_{n}(R^+))=F_*(W_{n-1}(R) \cap I_{n-1}^{[p^{e+1}]}W_{n-1}(R^+)),
\]
where the last equality follows from Lemma \ref{regular sequence case}~(1).
Thus we obtain
\[
[c'] F^e(\alpha) \in I_{n-1}^{[p^{e+1}]}W_{n-1}(R^+).
\]
Thus we obtain that
\[
[c^2c']F^e(\alpha) \in I_{n-1}^{[p^{e+1}]}
\]
by Lemma~\ref{regular sequence case}~(2) for every integer $e \geq 0$, where $c \in R^\circ$ is a test element.
In particular, we have $\alpha \in (I_{n-1}^{[p]})^*$.

In order to prove (3), let us show the following:
\begin{claim}\label{claim:colon-cap}
Let $c$ be a test element of $R$.
Then we have
\[
[c^2](I_n^{[p^e]} \colon [x_{r+1}^{p^{e}}]) \subseteq I_n^{[p^e]}
\]
for every integer $e \geq 0$.
\end{claim}
\begin{claimproof}
We take an integer $e \geq 0$.
Since $W_n(R^+)$ is a big Cohen-Macaulay $W_n(R)$-module by Proposition \ref{prop:big-CM}, we have
\[
(I_n^{[p^e]}W_n(R^+) \colon [x_{r+1}^{p^e}]) \subseteq I_n^{[p^e]}W_n(R^+).
\]
In particular, we have
\[
[c^2](I_n^{[p^e]} \colon [x_{r}^{p^e}]) \subseteq [c^2](W_n(R) \cap I_n^{[p^e]}W_n(R^+)) \subseteq I_n
\]
by Lemma~\ref{regular sequence case}~(2).
\end{claimproof}

Let us show (3). 
The first inclusion and the inclusion $I_n^* \subseteq (I_n^* \colon [x_{r+1}])$ are clear.
It is enough to show that $(I_n^* \colon [x_{r+1}]) \subseteq I_n^*$.
We take $\alpha \in (I_n^* \colon [x_{r+1}])$, then there exists $c_1 \in R^\circ$ such that we have
\[
[c_1]F^e([x_{r+1}]\alpha) \in I_n^{[p^e]}
\]
for every integer $e \geq 0$.
In particular, we have
\[
[c_1c^2]F^e(\alpha) \in [c^2](I_n^{[p^e]} \colon [x_{r+1}^{p^e}]) \overset{(\star)}{\subseteq} I_n^{[p^e]}, 
\]
where $(\star)$ follows from Claim \ref{claim:colon-cap}. 
Thus, we have $\alpha \in I_n^*$, as desired.
\end{proof}

\begin{proposition}[cf.~\cite{Smith94}*{Proposition~5.1}]\label{prop:tight-vs-plus}
We use \cref{n-local}.
Let $x_1,\ldots, x_r$ be a part of a system of parameters.
We consider an inductive system
\[
\left\{\frac{W_n(R)}{([x_1^m],\ldots,[x_d^m])}\right\}_{m \in \Z_{\geq 1}}
\]
with maps
\[
\frac{W_n(R)}{([x_1^m],\ldots,[x_d^m])} \to \frac{W_n(R)}{([x_1^{m+l}],\ldots,[x_d^{m+l}])}\ ;\ \overline{\alpha} \mapsto \overline{[x_1^l \cdots x_d^l]\alpha}.
\]
It is known that the colimit of the inductive system coincides with $H^d_{\m}(W_n(R))$.
Thus, we obtain the $W_n(R)$-module homomorphism 
\[
i^m_R \colon W_n(R)/([x_1^m],\ldots,[x_d^m]) \to H^d_{\m}(W_n(R)).
\]
We denote $i^1_R$ by $i_R$.
Let $n \geq 1$ be an integer and we set $I_n:=([x_1],\ldots, [x_r]) \subseteq W_n(R)$.
\begin{enumerate}
    \item We have $I_n^*/I_n=i_R^{-1}(\wt{0^*_n})$, where $\wt{0^*_n}$ is defined by \cite{KTTWYY3}*{Definition~3.16}.
    \item We have $I_n^*=I_n^+$.
\end{enumerate}
\end{proposition}

\begin{proof}
First, we prove (1).
We consider the following commutative diagram;
\begin{equation}\label{eq:cl:t-vs-p}
\begin{tikzcd}
    W_n(R)/I_n \arrow[r,"i_R"] \arrow[d,"F^e"] & H^d_{\m}(W_n(R)) \arrow[d,"F^e"] \\
    F^e_*W_n(R)/I_n^{[p^e]} \arrow[r,"i_R^{p^e}"] & F^e_*H^d_{\m}(W_n(R)).
\end{tikzcd}
\end{equation}
The inclusion $I_n^*/I_n \subseteq i_R^{-1}(\wt{0^*_n})$ follows from the diagram (\ref{eq:cl:t-vs-p}).
Furthermore, by the diagram, it is enough to show that there exists $c \in R^0$ such that $[c] \cdot \mathrm{Ker}(i_R^{p^e})=0$ for every $e \geq 0$.
The kernel of $i^{p^e}_R$ is
\[
\bigcup_{m \in \Z_{\geq 1}} ([x_1^{m+p^e}],\ldots,[x_d^{m+p^e}] \colon [x_1^m \cdots x_d^m]).
\]
Since we have
\begin{align*}
    ([x_1^{m+p^e}],\ldots,[x_d^{m+p^e}] \colon [x_1^m \cdots x_d^m]) 
    &\subseteq ([x_1^{m+p^e}],\ldots,[x_d^{m+p^e}] \colon [x_1^m \cdots x_d^m])^+ \\
    &\overset{(\star_1)}{=}(I^{[p^e]}_n)^+,
\end{align*}
where $(\star_1)$ follows from \cref{prop:big-CM}. 
By \cref{regular sequence case}~(2), there exists $c \in R^0$ such that $[c] (I^{[p^e]}_n)^+ \subseteq I_n^{[p^e]}$, as desired.

Next, to prove (2), we consider the map
\[
i_{R^+} \colon W_n(R^+)/I_nW_n(R^+) \to H^d_{\m}(W_n(R^+))
\]
constructed by a similar way to a construction of $i_R$.
Then we have a commutative diagram;
\begin{equation}
\begin{tikzcd}
    W_n(R)/I_n \arrow[r,"i_R"] \arrow[d,"a"] & H^d_{\m}(W_n(R)) \arrow[d,"b"] \\
%    W_n(R)/I_n^* \arrow[r,"\overline{i_R}"] \arrow[d,"a"] & H^d_{\m}(W_n(R))/\wt{0^*_n} \arrow[d,"b"] \\
    W_n(R^+)/I_nW_n(R^+) \arrow[r,"i_{R^+}"] & H^d_{\m}(W_n(R^+)).
\end{tikzcd}
\end{equation}
Since $i_{R^+}$ is injective, the kernel of $a$ coincides with $i_R^{-1}(\mathrm{Ker}(b))$.
By \cite{KTTWYY3}*{Theorem 5.6}, we have $\mathrm{Ker}(b)=\wt{0^{*}_n}$.
Therefore, we have
\[
I_n^+/I_n=\mathrm{Ker}(a)=i_R^{-1}(\wt{0^*_n}) \overset{(\star_3)}{=} I_n^*/I_n,
\]
where $(\star_3)$ follows from \cref{prop:tight-vs-plus}~(1).
Thus we have $I_n^+=I_n^*$, as desired.
\end{proof}

\begin{corollary}\label{cor:plus cl vs quasi-tight cl}
We use \cref{n-local}.
Let $x_1,\ldots, x_r$ be a part of a system of parameters.
Let $I_n:=([x_1],\ldots, [x_r]) \subseteq W_n(R)$  for integers $n \geq 1$.
Then we have
\[
I_n^*R={\rm Ker}(R \to Q^+_{R,n} \otimes W_n(R)/I_n),
\]
where $Q^+_{R,n}$ is defined by $Q^+_{R,1}:=R^+/R$ and
\[
Q^+_{R,n}:={\rm Coker}(F_*W_{n-1}(R) \xrightarrow{V} W_n(R) \to W_n(R^+))
\]
for $n \geq 2$.
\end{corollary}

\begin{proof}
By \cref{colon capturing} (2) and \cref{prop:tight-vs-plus}~(2), we obtain following commutative diagram in which each horizontal sequence is exact:
\[
\begin{tikzcd}
    0 \arrow[r]  & F_*W_{n-1}(R)/(I_{n-1}^{[p]})^* \arrow[r] \arrow[d,equal] & W_n(R)/I_n^* \arrow[r] \arrow[d,hookrightarrow] & R/I_n^*R \arrow[r]  & 0 \\
    0 \arrow[r] & F_*W_{n-1}(R)/(I_{n-1}^{[p]})^* \arrow[r] & W_n(R^+)/I_nW_n(R^+) \arrow[r] & Q^+_{R,n} \otimes W_n(R)/I_n \arrow[r] & 0.
\end{tikzcd}
\]
Thus, we obtain the injection $R/I_n^*R \hookrightarrow Q^+_{R,n} \otimes W_n(R)/I_n$, as desired.
\end{proof}

\subsection{Quasi-tightly closedness for systems of parameters and quasi-$F$-rationality}
In this subsection, we prove that quasi-$F$-rationality is equivalent to quasi-tightly closedness of parameter ideals (\cref{thm:new-old-h}), the result implies that quasi-$F$-rationality is a generalization of $F$-rationality.

First, we define the quasi-tightly closedness for a part of system of parameters.
\begin{definition}\label{defn:q-t.c}
We use \cref{n-local}.
Let $x_1,\ldots,x_r$ be a part of a system of parameters of $R$.
\begin{itemize}
    \item Let $h \geq 1$ be an integer. 
    We set $I_h^{\bm{v}}:=([x_1^{v_1}],\ldots,[x_r^{v_r}]) \subseteq W_h(R)$ for every integers $v_1,\ldots,v_r \geq 1$.
    We say that $x_1,\ldots,x_r$ is \emph{$h$-quasi-tightly closed} if for every integers $v_1,\ldots,v_r \geq 1$, then we have
    \[
    (I_h^{\bm{v}})^*R=I_1^{\bm{v}}.
    \]
    \item We say that $x_1,\ldots,x_r$ is \emph{quasi-tightly closed} if $x_1,\ldots,x_r$ is $h$-quasi-tightly closed for some positive integer $h$.
\end{itemize}
\end{definition}

\begin{remark}
In \cref{defn:q-t.c}, the ideal $(I_n^{\bm v})^*R$ depends on the choice of $x_1,\ldots,x_r$ (cf.~\cref{ex} (2)), thus we do not say "$I$ is quasi-tightly closed''.
\end{remark}

\begin{lemma}\label{lem:inf-limit-result}
We use \cref{n-local}.
Let $x_1,\ldots,x_r$ be a part of a system of parameters of $R$ and we set $I_n^{\bm{v}}:=([x_1^{v_1}],\ldots,[x_r^{v_r}]) \subseteq W_n(R)$ for an integer $n \geq 1$ and $\bm{v}=(v_1,\ldots,v_r) \in \Z_{\geq 1}^{r}$.
If $x_1,\ldots,x_r$ is $h$-quasi-tightly closed, then  we have
\begin{enumerate}
    \item $(I_{n+h-1}^{\bm{v}})^*W_n(R)=I_{n}^{\bm{v}}$ for integers $n \geq 1$ and
    \item $V^{-1}(I_n^{\bm{v}})=F_*I_{n-1}^{p\bm{v}}$ for integers $n \geq 2$.
\end{enumerate}
\end{lemma}

\begin{proof}
The assertion (1) for $n=1$ is clear.
Next, we assume $n \geq 2$ and take an integer $e \geq 0$.
Let us show that the assertion (1) for $n-1$ implies the assertions (1) and (2) for $n$.
We consider the following commutative diagram in which each horizontal sequence is exact:
\begin{equation}
\begin{tikzcd}
    0 \arrow[r] & V^{-1}(I_{n+h-1}^{\bm{v}}) \arrow[r,"V"] \arrow[d,"R^{h-1}"] & I_{n+h-1}^{\bm{v}} \arrow[r,"R^{n+h-2}"] \arrow[d,"R^{h-1}"] & I_1^{\bm{v}} \arrow[r] \arrow[d,equal] & 0 \\
    0 \arrow[r] & V^{-1}(I_n^{\bm{v}}) \arrow[r,"V"] & I_n^{\bm{v}} \arrow[r,"R^{n-1}"] & I_1^{\bm{v}} \arrow[r] & 0.
\end{tikzcd}
\end{equation}
By the snake lemma, we have
\begin{equation}\label{eq:V-inverse}
   V^{-1}(I_n^{\bm{v}})=R^{h-1}(V^{-1}(I_{n+h-1}^{\bm{v}})) \overset{(\star_1)}{\subseteq} F_*(I_{n+h-2}^{p\bm{v}})^*W_{n-1}(R) \overset{(\star_2)}{=}F_*I_{n-1}^{p\bm{v}}, 
\end{equation}
where $(\star_1)$ follows from Lemma \ref{colon capturing} (2) and $(\star_2)$ follows from (1) for $n-1$.
Since the converse inclusion is clear, we obtain (2) for $n$.
Next, we consider the following commutative diagram in which each horizontal sequence is exact:
\begin{equation}
    \begin{tikzcd}
        0 \arrow[r] & F_*(I_{n+h-2}^{p\bm{v}})^* \arrow[r,"V"] \arrow[d,"R^{h-1}"] & (I_{n+h-1}^{\bm{v}})^* \arrow[r,"R^{n+h-2}"] \arrow[d,"R^{h-1}"] & (I_{n+h-1}^{\bm{v}})^*R \arrow[r] \arrow[d] & 0 \\
        0 \arrow[r] & F_*W_{n-1}(R) \arrow[r,"V"] & W_n(R) \arrow[r,"R^{n-1}"] & R \arrow[r] & 0.
    \end{tikzcd}
\end{equation}
By the snake lemma, we obtain the exact sequence
\begin{equation}\label{eq:exact-m}
    0 \to F_*W_{n-1}(R)/(I_{n+h-2}^{p\bm{v}})^*W_{n-1}(R) \xrightarrow{V} W_n(R)/(I_{n+h-1}^{\bm{v}})^*W_n(R) \xrightarrow{R^{n-1}} R/(I_{n+h-1}^{\bm{v}})^*R \to 0. 
\end{equation}
We take $\alpha \in (I_{n+h-1}^{\bm{v}})^*W_n(R)$, then we have
\[
R^{n-1}(\alpha) \in  (I_{n+h-1}^{\bm{v}})^*R=I_1^{\bm{v}}
\]
by the induction hypothesis.
Then there exists $\alpha' \in I_n^{\bm{v}}$ such that $R^{n-1}(\alpha-\alpha')=0$.
In particular, there exists $\beta \in W_{n-1}(R)$ such that $V\beta=\alpha-\alpha'$.
We have
\[
\beta \in V^{-1}( (I_{n+h-1}^{\bm{v}})^*W_n(R)) \overset{(\star_3)}{=} F_*(I_{n+h-2}^{p\bm{v}})^*W_{n-1}(R) \overset{(\star_4)}{=} F_*I_{n-1}^{p\bm{v}},
\]
where $(\star_3)$ follows from (\ref{eq:exact-m}) and $(\star_4)$ follows from (1) for $n-1$.
Therefore, we have $\alpha-\alpha' \in I_n^{\bm{v}}$.
Since $\alpha'$ is contained in $I_n^{\bm{v}}$, so is $\alpha$, thus we obtain $I_n^{\bm{v}} =(I_{n+h-1}^{\bm{v}})^*W_n(R)$, as desired.
\end{proof}

\begin{lemma}\label{lem:descent}
We use \cref{n-local}.
Let $x_1,\ldots,x_r$ be a part of a system of parameters of $R$ and we set  $J_n:=([x_1],\ldots,[x_{r-1}]) \subseteq W_n(R)$ for integers $n \geq 1$.
If $x_1,\ldots,x_r$ is $h$-quasi-tightly closed, then $J_{2h-1}^*R=J_1$ and $(J_1 \colon x_r)=J_1$.
\end{lemma}

\begin{proof}
First, let us show the following claim.
\begin{claim}\label{claim:desent-m}
We have $J_{2h-1}^*R \subseteq x_r^m J_h^*R+J_1$ for every positive integer $m$.
\end{claim}

\begin{proof}[Proof of Claim]
We set $I_{n,m}:=(J_n,[x_r^m]) \subseteq W_n(R)$.
Since we have
\[
I_{n,m}^{[p^e]} \subseteq ((J_{n+h-1,m}^{[p^e]})^*W_n(R),[x_r]^{mp^e}) \subseteq (I_{n+h-1,m}^{[p^e]})^*W_n(R) \overset{(\star_1)}{=} I_{n,m}^{[p^e]},
\]
where $(\star_1)$ follows from Lemma \ref{lem:inf-limit-result}, we obtain
\begin{equation}\label{eq:gene}
((J_{n+h-1}^{[p^e]})^*W_n(R),[x_r]^{mp^e})=I_{n,m}^{[p^e]}  
\end{equation}
for integers $n \geq 1$ and $e \geq 0$.
We consider the following commutative diagram in which each horizontal sequence is exact:
\begin{equation*} 
    \begin{tikzcd}[column sep=0.5cm]
    0 \arrow{r} & F_*W_{h-1}(R)/(J_{2h-2}^{[p]})^*W_{h-1}(R) \arrow{r}{V} \arrow{d}{\cdot [x_r]^{mp}} & W_h(R)/J_{2h-1}^*W_h(R) \arrow{r}{R^{h-1}} \arrow{d}{\cdot [x_r]^m} & R/J_{2h-1}^*R \arrow{r} \arrow{d}{\cdot x_r^m} & 0 \\
    0 \arrow{r} & F_*W_{h-1}(R)/(J_{2h-2}^{[p]})^*W_{h-1}(R) \arrow{r}{V} & W_h(R)/J_{2h-1}^*W_h(R) \arrow{r}{R^{h-1}} & R/J_{2h-1}^*R  \arrow{r} & 0.
    \end{tikzcd}
\end{equation*}
By (\ref{eq:gene}), the cokernel of the middle vertical map is
\[
W_h(R)/(J_{2h-1}^*W_n(R),[x_r]^m)=W_h(R)/I_{h,m}
\]
and the cokernel of the left vertical map is
\[
F_*W_{h-1}(R)/((J_{2h-2}^{[p]})^*W_{h-1}(R),[x_r]^{mp})=F_*W_{h-1}(R)/I_{h-1,m}^{[p]}.
\]
By the snake lemma and Lemma \ref{lem:inf-limit-result}, we have 
\begin{equation}\label{eq:key-surj}
   (J_{2h-1}^*W_{h}(R) \colon [x_r]^m)R= (J_{2h-1}^*R \colon x_r^m). 
\end{equation}
We take $a \in J_{2h-1}^*R$.
Since we have
\[
J_{2h-1}^*R \subseteq I_{h,m}^*R=I_{1,m}=(J_1,x_r^m),
\]
there exist $j \in J_1$ and $b \in R$ such that $a=x_r^mb+j$.
In particular, we have
\[
b \in (J_{2h-1}^*R \colon x_r^m) \overset{(\star_2)}{=} (J_{2h-1}^*W_h(R) \colon [x_r]^m)R \subseteq (J_h^* \colon [x_r]^m)R \overset{(\star_3)}{=} J_{h}^*R,  
\]
where $(\star_2)$ follows from (\ref{eq:key-surj}) and $(\star_3)$ follows from Lemma \ref{colon capturing} (3).
Therefore, we obtain $J_{2h-1}^*R \subseteq x_r^m J_{h}^*R+J_1$.
\end{proof}
\noindent
By Claim \ref{claim:desent-m}, for every positive integer $m$, we have $J_{2h-1}^*R/J_1 \subseteq x_r^mJ_h^*R/J_1$.
Since $x_r$ is contained in the maximal ideal, we have $J_{2h-1}^*R=J_1$.

Next, let us show $(J_1 \colon x_r)=J_1 $.
It is enough to show that $(J_1 \colon x_r) \subseteq J_{2h-1}^*R$.
% We have
% \begin{align*}
%     (J_1 \colon x_r) =(J_{2h-1}^*R \colon x_r) \overset{(\star_4)}{=} (J_{2h-1}^*W_h(R) \colon [x_r])R \subseteq (J_{2h-1}^* \colon [x_r])R \overset{(\star_5)}{=}J_{2h-1}^*R=J_1,
% \end{align*}
% where $(\star_4)$ follows from (\ref{eq:key-surj}) and $(\star_5)$ follows from \ref{colon capturing(3)}, as desired.
If $h=1$, it follows from Lemma~\ref{colon capturing}~(3).
We assume $h \geq 2$.
We consider the following commutative diagram in which each horizontal sequence is exact:
\begin{equation}\label{eq:exact-J-I}
    \begin{tikzcd}
    0 \arrow{r} & F_*W_{2h-2}(R)/V^{-1}(J_{2h-1}) \arrow{r}{V} \arrow{d}{\cdot [x_r]^p} & W_{2h-1}(R)/J_{2h-1} \arrow{r}{R^{2h-2}} \arrow{d}{\cdot [x_r]} & R/J_1 \arrow{r} \arrow{d}{\cdot x_r} & 0 \\
    0 \arrow{r} & F_*W_{2h-2}(R)/V^{-1}(J_{2h-1}) \arrow{r}{V} & W_{2h-1}(R)/J_{2h-1} \arrow{r}{R^{2h-2}} & R/J_1  \arrow{r} & 0.
    \end{tikzcd}
\end{equation}
Furthermore, we have
\[
V^{-1}(I_{2h-1,1}) \overset{(\star_4)}{=} F_*I_{2h-2,1}^{[p]} \subseteq (V^{-1}(J_{2h-1}),[x_r]) \subseteq V^{-1}(I_{2h-1,1})
\]
where $(\star_4)$ follows from Lemma \ref{lem:inf-limit-result}, thus the cokernel of the left vertical map in (\ref{eq:exact-J-I}) is
\[
F_*W_{2h-2}(R)/V^{-1}(I_{2h-1,1}).
\]
Moreover, since we have $(J_{2h-1},[x_r])=I_{2h-1,  1}$, applying the snake lemma for (\ref{eq:exact-J-I}), we obtain 
\begin{equation}\label{eq:des-1}
    (J_{2h-1} \colon [x_r])R= (J_1 \colon x_r).
\end{equation}
By Lemma \ref{colon capturing} (3), we have
\begin{equation}\label{eq:des-2}
    (J_{2h-1} \colon [x_r])R \subseteq J_{2h-1}^*R.
\end{equation}
Combining (\ref{eq:des-1}) and (\ref{eq:des-2}), we obtain $(J_1 \colon x_r) \subseteq J_{2h-1}^*R$,
as desired.
\end{proof}

\begin{theorem}\label{thm:new-old-inf}
We use \cref{n-local}.
Let $x_1,\ldots,x_d$ be a system of parameters of $R$.
Then the following are equivalent to each other.
\begin{enumerate}
    \item there exists a positive integer $h$ such that $x_1,\ldots,x_d$ is $h$-quasi-tightly closed, and
    \item there exists a positive integer $h'$ such that  $x_1,\ldots,x_r$  is  $h'$-quasi-tightly closed for every $1 \leq r \leq d$.
\end{enumerate}
\end{theorem}

\begin{proof}
The assertion (2) $\Rightarrow$ (1) is obvious.
Let us show (1) $\Rightarrow$ (2).
For an integer $1 \leq r \leq d$, we set $h_r:=2^{d-r}(h-1)+1$.
We consider the following assertion for integers $1 \leq r \leq d$:
\begin{enumerate}
    \item[($\star_r$)]  $x_1,\ldots,x_r$ is $h_r$-quasi-tightly closed.
\end{enumerate}
Since we have $h_d=h-1+1=h$, the assertion($\star_d$) holds by (2).
Thus, it is enough to show that  $(\star_r)$ implies $(\star_{r-1})$ for $2 \leq r \leq d$.
We take an integer $2 \leq r \leq d$ and assume $(\star_r)$.
We set $J_n^{\bm{v}}:=([x_1]^{v_1},\ldots,[x_{r-1}]^{{v_{r-1}}}) \subseteq W_n(R)$ and $I_n^{\bm{v}}:=([x_1]^{v_1},\ldots,[x_r]^{v_r}) \subseteq W_n(R)$ for a positive integer $n$ and $\bm{v}:=(v_1,\ldots,v_r) \in \Z_{\geq 1}^{r}$.
Since $x_1^{v_1},\ldots,x_r^{v_r}$ is $h_r$-quasi-tightly closed, by Lemma \ref{lem:descent}, we have $(J_{2h_{r}-1}^{\bm{v}})^*R=J_1^{\bm{v}}$ for every $\bm{v} \in \Z_{\geq 1}^r$.
Since $2h_{r}-1=2^{d-(r-1)}(h-1)+2-1=h_{r-1}+1$, the assertion ($\star_{r-1}$) holds.
\end{proof}

\begin{theorem}[Theorem~\ref{intro:thm:new-old-h}]\label{thm:new-old-h}
We use \cref{n-local}.
Let $h \in \Z_{>0}$.
Then the following are equivalent to each other.
\begin{enumerate}
    \item  $R$ is $h$-quasi-$F$-rational (cf.~\cite{KTTWYY3}*{Definition~2.5}), 
    \item for every system of parameters $x_1,\ldots,x_d$ is $h$-quasi-tightly closed.
    \item there exists a system of parameters $x_1,\ldots,x_d$ of $R$ such that $x_1,\ldots,x_d$ is $h$-quasi-tightly closed.
\end{enumerate}
\end{theorem}

\begin{proof}
First, let us show (1) $\Rightarrow$ (2).
We take a system of parameters $x_1,\ldots,x_d$, and let $I_h:=([x_1],\ldots,[x_d]) \subseteq W_n(R)$ for every positive integer $h$.
Since $R$ is Cohen-Macaulay, the map
\[
i_R \colon W_h(R)/I_h \to H^d_\m(W_h(R))
\]
taking colimit is injective (cf.~proof of \cref{prop:tight-vs-plus}).
Moreover, we have $i_R^{-1}(0^*_h)=I_h^*/I_h$ by the proof of \cref{prop:tight-vs-plus}.
Since $R$ is $h$-quasi-$F$-rational, we have $R^{h-1}(0^*_h) =0$.
Therefore, we have
\[
I_h^*R=I_1+i_R^{-1}(0^*_h)R=I_1,
\]
as desired.
The implication (2) $\Rightarrow$ (3) is clear.

Finally, let us show (3) $\Rightarrow$ (1).
By Theorem \ref{thm:new-old-inf}, for every $1 \leq r \leq d$, $x_1,\ldots,x_r$ is $h'$-quasi-tightly closed for some $h'$.
Let us show that $x_1,\ldots,x_r$ is a regular sequence for every $1 \leq r \leq d$ by induction on $r$.
Since $R$ is a domain, $x_1$ is a regular element.
We assume $r \geq 2$ and $x_1,\ldots,x_{r-1}$ is a regular sequence. 
By Lemma \ref{lem:descent}, we have $((x_1,\ldots,x_{r-1}) \colon x_r)=(x_1,\ldots,x_{r-1})$, thus $x_1,\ldots,x_{r}$ is a regular sequence.
Therefore, we obtain that $x_1,\ldots,x_d$ is a regular sequence, and in particular, $R$ is Cohen-Macaulay.
It is enough to show that $0^*_{h}=0$.
We note that $H^d_\m(W_h(R))$ is the colimit of $W_h(R)/([x_1]^m,\ldots,[x_d]^m)$.
Therefore, $\wt{0^*_h}$ is also the colimit of 
\[
([x_1]^m,\ldots,[x_d]^m)^*/([x_1]^m,\ldots,[x_d]^m).
\]
Thus, $0^*_h$ is the colimit of $([x_1]^m,\ldots,[x_d]^m)^*R/(x_1^m,\ldots,x_d^m)=0$.
\end{proof}

\begin{corollary}\label{cor:power-t-c}
We use \cref{n-local}.
We assume that $R$ is $n$-quasi-$F$-rational for some positive integer $n$.
Let $x_1,\ldots,x_r$ be a part of a system of parameters and $I_n:=([x_1],\ldots,[x_r]) \subseteq W_n(R)$.
Let $\lambda \geq 1$ be an integer.
Then $I_1^\lambda=(I_n^{\lambda})^*R$.
\end{corollary}

\begin{proof}
Let $\alpha \in (I_n^{\lambda})^*$.
Then there exists $c \in R^\circ$ such that $[c]F^e(\alpha) \in (I_n^\lambda)^{[p^e]}$ for every $e \geq 1$.
By \cite{LT}*{Section 3}, we have
\[
[c]F^e(\alpha) \in (I_n^\lambda)^{[p^e]}=\bigcap_{\lambda_1+\cdots+\lambda_r=\lambda+r-1} ([x_1]^{\lambda_1},\ldots,[x_r]^{\lambda_r})^{[p^e]},
\]
thus we have $\alpha \in ([x_1]^{\lambda_1},\ldots,[x_r]^{\lambda_r})^*$ for every $\lambda_1,\ldots,\lambda_r$ with $\lambda_1+\ldots+\lambda_r=\lambda+r-1$.
In particular, we have
\[
R^{n-1}(\alpha)=\bigcap_{\lambda_1+\cdots+\lambda_r=\lambda+r-1} ([x_1]^{\lambda_1},\ldots,[x_r]^{\lambda_r})^*R \overset{(\star)}{=} \bigcap_{\lambda_1+\cdots+\lambda_r=\lambda+r-1} (x_1^{\lambda_1},\ldots,x_r^{\lambda_r})=I_1^\lambda,
\]
where $(\star)$ follows from the assumption that $R$ is $n$-quasi-$F$-rational and \cref{thm:new-old-h}.
\end{proof}

\begin{corollary}\label{t.c on witt}
We use  \cref{n-local}.
Let $n$ be a positive integer.
We assume that for every parameter ideal $I$ on $W_n(R)$, we have $I^*=I$.
Then $R$ is $F$-rational.
\end{corollary}

\begin{proof}
Let $x_1,\ldots,x_d$ be a system of parameter of $R$ and $I_n:=([x_1],\ldots,[x_d]) \subseteq W_n(R)$.
Since $I_n^*=I_n$, we have $I_n^*R=(x_1,\ldots,x_d)$.
By \cref{thm:new-old-h}, $R$ is quasi-$F$-rational, and in particular, $x_1,\ldots,x_d$ is a regular sequence.
Therefore, we obtain
\[
(x_1^{p^e},\ldots,x_d^{p^e})^*\overset{(\star_1)}{=} V^{-(n-1)}(I_n^*)=V^{-(n-1)}(I_n) \overset{(\star_2)}{=} (x_1^{p^e},\ldots,x_d^{p^e}),
\]
where $(\star_1)$ follows from \cref{colon capturing} (2) and $(\star_2)$ follows from the fact that $x_1,\ldots,x_d$ is a regular sequence.
\end{proof}

\begin{remark}
By Corollary~\ref{t.c on witt}, the tightly closedness of all parameter ideals in $W_h(R)$ is strictly stronger than $h$-quasi-$F$-rationality, which is equivalent to the $h$-quasi-tightly closedness of all systems of parameters.
\end{remark}

\subsection{Relation between Quasi-tightly closedness and quasi-$F$-regularity}
In this subsection, we study quasi-tight closure of ideals on quasi-$F$-regular ring and we use notations in \cite{KTTWYY3}*{Section~2.4}

\begin{proposition}\label{prop:gene t-c vs i-c}
Let $R$ be an $F$-finite normal domain in characteristic $p >0$ and $n$ a positive integer.
Let $I$ be an ideal of $R$ and $I_n$ an ideal of $W_n(R)$ such that $I_nR=I$.
Then we have
\[
T_n(\tau_n(R)I_n^*R) \subseteq I_nW_n\omega_R(-K_R),
\]
where $T_n$ is a map $T_n \colon R \simeq \omega_R(-K_R) \to W_n\omega_R(-K_R)$ induced by the trace map and $\tau_n(R)$ is defined in \cite{KTTWYY3}*{Definition~4.23}.
\end{proposition}

\begin{proof}
We take $a \in \tau_n(R)$ and $\alpha \in I_n^*$, then there exists $c \in R^\circ$ such that $[c]F^e(\alpha) \in I_n^{[p^e]}$ for every $e \geq 0$.
By \cite{KTTWYY3}*{Proposition~4.25}, there exists $c' \in R^\circ$ and $\omega \in W_n\omega_R(-K_R)$ such that $T^{e,cc'}_n(F^e_*\omega)=T_n(a)$, where $T^{e}_n$ is the $W_n\omega_R$-dual of $F^e \colon W_n(R) \to W_n(R)$ and
\[
T^{e,x}_n(F^e_*\omega)=T^e_n(F^e_*([x]\omega))
\]
for $x \in R$.
Therefore, we have
\begin{align*}
    T^{e,cc'}_n(F^e_*(F^e(\alpha) \omega)) &= T^{e}_n(F^e_*([cc']F^e(\alpha) \omega)) \subseteq T^e_n(F^e_*(I_n^{[p^e]}W_n\omega_R(-K_R))) \subseteq I_nW_n\omega_R(-K_R) \\
    T^{e,cc'}_n(F^e_*(F^e(\alpha) \omega)) &=\alpha T^{e,cc'}_n(F^e_*\omega) = \alpha T_n(a) = T_n(R^{n-1}(\alpha)a).
\end{align*}
Thus, we have $T_n(R^{n-1}(\alpha)a) \in I_nW_n\omega_R(-K_R)$, as desired.
\end{proof}

If $R$ is $n$-quasi-$F$-regular for some $n$, then we have $T_n(I_n^*R) \subseteq I_nW_n\omega_R$ by Proposition \ref{prop:gene t-c vs i-c}.
In particular, we have
\[
I_n^*/I_n \subseteq \Ker(R/I \simeq \omega_R(-K_R) \otimes W_n(R)/I_n \to W_n\omega_R(-K_R) \otimes W_n(R)/I_n).
\]
Therefore, the obstruction of quasi-tightly closedness of ideals on a quasi-$F$-regular ring is the kernel of the map $\omega_R(-K_R) \otimes W_n(R)/I_n \to W_n\omega_R(-K_R) \otimes W_n(R)/I_n$.
By \cref{ex} (3), it is not injective in general lift $I_n$.
Thus, the following question is natural to consider.

\begin{question}\label{ques:q-t.c on qFr}
Let $R$ be an $F$-finite normal domain in characteristic $p >0$ and $n$ a positive integer.
Let $I$ be an ideal of $R$.
Is there a lift $I_n$ on $W_n(R)$ such that the map $\omega_R(-K_R) \otimes W_n(R)/I_n \to W_n\omega_R(-K_R) \otimes W_n(R)/I_n$ is injective?
\end{question}

\subsection{Quasi-tight closure and integrally closure}
In this subsection, we study relationship between quasi-tight closure and integrally closure of ideals.
The results of this subsection were inspired by a discussion with Jakub Witaszek.

\begin{lemma}\label{lem:ideal-t-l}
Let $R$ be an $F$-finite Noetherian ring of characteristic $p >0$.
Let $f_1,\ldots,f_r \in R$ and $I_n:=([f_1],\ldots,[f_r]) \subseteq W_n(R)$ for positive integer $n$.
Let $n,e,\lambda$ be integers with $e \geq 0$ and $n,\lambda \geq 1$.
If $(g_0,g_1,\ldots,g_{n-1}) \in W_n(R)$ satisfies $g_i \in I_1^{(\lambda+r-1)p^{e+i}}$ for every $n-1 \geq i \geq 0$, then $(g_0,g_1,\ldots,g_{n-1}) \in (I_n^{\lambda})^{[p^{e}]}$.
\end{lemma}

\begin{proof}
We prove the assertion by induction on $n$.
In the case of $n=1$, it is clear.
We assume $n \geq 2$.
Let $\lambda':=\lambda+r-1$.
Since $g_0 \in I_1^{\lambda' p^{e}}$, we have a decomposition $g_0=\sum a_{\bm{k}}f_1^{k_1} \cdots f_r^{k_r}$, where the sum runs over $\bm{k}=(k_1,\ldots,k_r) \in \Z^r_{\geq 0}$ with $k_1+\cdots+k_r=\lambda'p^e$.
Then we have
\[
[g_0]=\left[\sum_{k_1+\cdots+k_r=\lambda'p^e}a_{\bm{k}}f_1^{k_1}\cdots f_r^{k_r}\right]=\sum_{k_1+\cdots+k_r=\lambda'p^e}[a_{\bm{k}}f_1^{k_1}\cdots f_r^{k_r}] +(0,g_1',\ldots,g'_{n-1}).
\]
By the proof of \cite{KTTWYY2}*{Proposition 7.5}, $g_i'$ is a sum of $p^i$-products of  elements in $\{a_{\bm{k}} f_1^{k_1}\cdots f_r^{k_r} \mid k_1+\ldots+k_r=\lambda'p^e \}$.
Therefore, we have $g_i' \in I^{rp^{e+i}}_1$.
On the other hand, we have
\[
(g_0,\ldots,g_{n-1})=\sum_{k_1+\cdots+k_r=\lambda'p^e}[a_{\bm{k}}f_1^{k_1}\cdots f_r^{k_r}]+V(g_1+g_1',\ldots,g_{n-1}+g_{n-1}').
\]
Since $g_i+g_i' \in I_1^{\lambda'p^{e+i}}$, by the induction hypothesis, we have
\[
(g_1+g_1',\ldots,g_{n-1}+g_{n-1}') \in (I_{n-1}^{\lambda})^{[p^{e+1}]}.
\]
For $k_1+\ldots+k_r=\lambda'p^e$, we have $[a_{\bm{k}}f_1^{k_1}\cdots f_r^{k_r}] \in (I_n^{\lambda})^{[p^{e}]}$, thus we have $(g_0,\ldots,g_{n-1}) \in (I_n^{\lambda})^{[p^e]}$, as desired. 
\end{proof}

\begin{theorem}[Theorem~\ref{thm:intro-gene t-c vs i-c}]\label{prop:int-cl vs tight-cl}
Let $R$ be an $F$-finite Noetherian ring of characteristic $p >0$.
Let $f_1,\ldots,f_r \in R$, $I_n:=([f_1],\ldots,[f_r]) \subseteq W_n(R)$ and $I:=(f_1,\ldots,f_r)$.
Then we have $\overline{I^{(\lambda+r-1)}} \subseteq (I_n^{\lambda})^*R$ for every integer $ \lambda \geq 1$.
\end{theorem}

\begin{proof}
Take $a \in \overline{I^{(\lambda+r-1)}}$.
Then there exists non-zero divisor $c \in R$ such that $ca^m \in I^{(\lambda+r-1)m}$ for every $m$.
By \cref{lem:ideal-t-l}, we have $[ca^{p^e}] \in (I_n^{\lambda})^{[p^e]}$ for every integers $n \geq 1$ and $e \geq 0$.
Therefore, we obtain that $[a]$ is contained in $(I_n^{\lambda})^*$, and in particular, $a$ is contained in $(I_n^{\lambda})^*R$, as desired.
\end{proof}

\begin{remark}
By \cref{prop:int-cl vs tight-cl} and \cref{cor:power-t-c}, we obtain that $\overline{I^{(\lambda+r-1)}} \subseteq I^\lambda$, where $I$ is generated by $r$ elements which is a part of a system of parameter of a quasi-$F$-rational ring.
However, the result was obtained by \cite{LT}*{Corollary~2.2} and \cite{KTTWYY3}*{Theorem~3.44}.
\end{remark}

% \begin{theorem}\label{thm:Brianson-Skoda}
% We work in the general setting (\cref{n-local}).
% We assume that $R$ is quasi-$F$-rational.
% Let $x_1,\ldots,x_r$ be a part of a system of parameters and $I:=(x_1,\ldots,x_r)$.
% Then $\overline{I^{(\lambda+r-1)}} \subseteq I^\lambda$.
% \end{theorem}

% \begin{proof}
% We take a positive integer $n$ such that $R$ is $n$-quasi-$F$-rational.
% By \cref{prop:int-cl vs tight-cl}, we have $\overline{I^{(\lambda+r-1)}} \subseteq (I_n^{\lambda})^*R$.
% Since $R$ is quasi-$F$-rational, by \cref{cor:power-t-c}, we have $(I_n^{\lambda})^*R=I^\lambda$.
% Thus, we have $\overline{I^(\lambda+r-1)} \subseteq I^\lambda$, as desired.
% \end{proof}

\section{Examples}
In this section, we introduce examples related to the definition of an analog of tight closure via Witt rings.
\cref{ex} (2) means that ''quasi-tight closure" used in the paper depends on the choice of generators.
By Proposition \ref{prop:gene t-c vs i-c}, the obstruction of ''quasi-tight closedness" on a quasi-$F$-regular ring is the kernel of the map $\omega_R(-K_R) \otimes W_n(R)/I_n \to W_n\omega_R(-K_R) \otimes W_n(R)/I_n$. 
\cref{ex} (3) implies that the map is not injective in general.

\begin{example}\label{ex}
Let $k$ be an algebraically closed field of characteristic $2$.
\begin{enumerate}
    \item Let $R:=k[x,y,z]/(x^3+y^3+z^3)$ and $I=(x,y)$.
    We take elements $f,g \in R$ such that $I=(f,g)$ and set $I_n:=([f],[g])$.
    Then for every positive integer $n$, we have $I_n^*R=(x,y,z^2)$.
    Indeed, we see that $I^*=(x,y,z^2)$, thus we have $I_n^*R \subseteq I^*=(x,y,z^2)$.
    On the other hand, we have
    \[
    (x,y,z^2) \subseteq \overline{I^2} \overset{(\star)}{\subseteq} I_n^*R,
    \]
    where $(\star)$ follows from \cref{prop:int-cl vs tight-cl}.
    \item Let $R:=k[x,y,z]/(z^4+x^5+y^8)$ and $I:=(x,y)$.
    We consider the two ideals $I_2:=([x],[y])$ and $I_2':=([x+y],[y])$.
    Then we have $z \in I_2^*R$ and $z \notin (I_2')^*R$, and in particular, $(I_2)^*R$ depends on the choice of the lifts of generators.
    Indeed, in $R^+$, we have $z=x^{5/4}+y^{2}$.
    Thus, we have
    \begin{align*}
        [z]&=[x^{5/4}]+[y^2]+(0,x^{5/4}y^2) \in I^*_2=I_2W_2(R^+) \cap W_2(R) \\
        [z] &=[x+y][x^{1/4}]+[y][x^{1/4}+y]+(0,(x+y)x^{1/4}y(x^{1/4}+y)) \\
        &\equiv (0,x^{3/2}y) \mod I_2'W_2(R^+).
    \end{align*}
    Thus, we have $z \in I_2^*R$.
    On the other hand, if $z$ is contained in $(I_2')^*R$, then there exists $a \in R$ such that $a+x^{3/2}y$ is contained in $I^{[2]}R^+$, and in particular, $a^2+x^3y^2$ is contained in $I^{[4]}R^+ \cap R=(I^{[4]})^*$.
    However, since we see that it is not true, thus we have $z \notin (I_2')^*R$.
    \item Let $R:=k[x,y]_{(x,y)}$ and $I_2:=([x],[y],[x+y]) \subseteq W_2(R)$.
    Then $\omega_R \otimes W_2(R)/I_2 \to W_2\omega_R \otimes W_2(R)/I_2$ induced by the trace map is not injective.
    Indeed, the dual of the map is $\pi \colon H^3_{(x,y)}(W_2(R))[I_2] \to H^3_{(x,y)}(R)[I_2]$, where $(-)[I_2]$ means the submodule of $I_2$-torsion elements.
    It is enough to show that $\pi$ is not surjective.
    In particular, it is enough to show that every lift of $[1/xy]$ is not $I_2$-torsion.
    Suppose there exists a lift $\alpha \in H^3_{(x,y)}(W_2(R))[I_2]$.
    Since $\alpha$ is $([x],[y])$-torsion, there exists $a \in R$ such that $\alpha=[(1,a)/[xy]]$.
    Since we have
    \[
    \alpha [x+y]=\left[\frac{(1,a)}{[xy]}\right][x+y]=\left[\frac{(0,xy)}{[xy]}\right] \neq 0,
    \]
    thus $\alpha$ is not $I_2$-torsion.
    Therefore, $\omega_R \otimes W_2(R)/I_2 \to W_2\omega_R \otimes W_2(R)/I_2$ is not injective.
\end{enumerate}
\end{example}

\bibliographystyle{skalpha}
\bibliography{bibliography.bib}

\end{document}